\documentclass[a4]{amsart}
\usepackage{cite}
\usepackage{amssymb,amscd,amsmath,enumerate,verbatim,calc}
\usepackage{mathtools}
\usepackage{mathrsfs}
\numberwithin{equation}{section}

\allowdisplaybreaks[4]

\DeclarePairedDelimiter\floor{\lfloor}{\rfloor}

\usepackage[a4paper, total={6in, 8in}]{geometry}

\newtheorem{lemma}{Lemma}[section]
\newtheorem{theorem}{Theorem}[section]
\newtheorem{definition}{Definition}[section]
\newtheorem{proposition}{Proposition}[section]
\newtheorem{solution}{Solution}[section]
\newtheorem{example}{Example}[section]
\theoremstyle{remark}

\newcommand{\Z}{\ensuremath{\mathbb{Z}}}

\newcommand{\R}{\ensuremath{\mathbb{R}}}

\newcommand{\N}{\ensuremath{\mathbb{N}}}

\begin{document}
	
	\title[]{Constructing Multiresolution Analysis via Wavelet Packets on Sobolev Space in Local Fields}
	
	\author{Manish Kumar}
	
	\address{(M. Kumar) Department of Mathematics\\ Birla Institute of Technology and Science-Pilani, Hyderabad\\ Telangana 500078, India}
	
	\email{manish.math.bhu@gmail.com; manishkumar@hyderabad.bits-pilani.ac.in}

	
	
	
	\date{\today}
	
	
	
	\begin{abstract}

We define Sobolev spaces $H^{\mathfrak{s}}(K_q)$ over a local field $K_q$ of finite characteristic $p>0$, where $q=p^c$ for a prime $p$ and $c\in \mathbb{N}$. This paper introduces novel fractal functions, such as the Weierstrass type and 3-adic Cantor type, as intriguing examples within these spaces and a few others. Employing prime elements, we develop a Multi-Resolution Analysis (MRA) and examine wavelet expansions, focusing on the orthogonality of both basic and fractal wavelet packets at various scales. We utilize convolution theory to construct Haar wavelet packets and demonstrate the orthogonality of all discussed wavelet packets within $H^{\mathfrak{s}}(K_q)$, enhancing the analytical capabilities of these Sobolev spaces.

	 	\end{abstract}
	
	\maketitle
	
	\section{Introduction}
  Wavelet packets (sometimes known as wave packets or subband trees) are a further decomposition of wavelet space; see Coifman, Meyer, Wickerhauser, and others in \cite{coifman1990wavelet,coifman1992size,coifman1994signal,WICKERHAUSER1992679} for a background. Wavelet packets provide very elegant and sophisticated mathematical tools. They are essential for applications in science and engineering while dealing with both the development of theory and applications. It was already used in signal processing and compression technique; see \cite{WICKERHAUSER1992679,coifman1994signal}. Medical science uses it as a preclinical diagnosis (see \cite{zhang2015preclinical}). In \cite{coifman1994signal}, the authors investigated in-depth the theory of wavelet packets and their properties in a classical sense. A detailed analysis of the theory of wavelet packets in $L^{2}({\mathbb R})$ is given in the famous books by Chui \cite{chui1993wavelets}, Hern$\acute{a}$ndez and Weiss \cite{hernandez1996first}, or by Ruch and Fleet \cite{van2011wavelet}. Similarly, motivated by the work of Bastin and Laubin \cite{bastin1997regular}, Pathak and Kumar in \cite{pathak2015wavelet} constructed a MRA and derived many important properties of wavelet packets for the Sobolev space $H^{\mathfrak s}({\mathbb R})$. 
	The depth analysis of number theory on local fields can be seen in \cite{gouvea1997p, koblitz2012p, robert2000course, serre2013local, PAGANO2022398}. Several mathematicians focused on working with local fields and have found countless applications in many research areas, including number theory, representation theory, class field theory, division algebras, dynamical systems, quadratic forms, linear algebraic groups, and algebraic geometry. Also, many authors have investigated the class field theory (which is a fundamental branch of algebraic number theory) using Fourier analysis on the local field, and certain zeta functions that can be found in \cite{connes1997trace, connes1999trace, ramakrishnan1998fourier, tate1997fourier}. To the best of my knowledge, a minimal amount of work has been researched for wavelet packets in classical spaces, say space of all square-integrable functions associated with a local field $K_q$ of finite characteristic $p>0$, denoted by $L^{2}({K_q})$, see \cite{behera2012wavelet,pathak2023existence}, but not for distributional spaces, for instance, {\em{Sobolev space associated with local field $K_q$ of finite characteristic $p>0$ denoted by $H^{\mathfrak s}(K_q)$}}. Moreover, reference \cite{behera2012wavelet} discusses only theoretical aspects without providing examples, and examples of wavelet packets in Sobolev space $H^{\mathfrak s}(K_q)$ are absent from the current literature.\\
\par	Motivated by this fact, including the lack of examples, we wish to define Sobolev space using $K_q$ and illustrate it via a few interesting examples, including fractals. Further, we construct MRA on Sobolev space $H^{\mathfrak s}(K_q)$. Following the technique of MRA, we derive wavelet packets and investigate their important properties with exciting new examples in $H^{\mathfrak s}(K_q)$. Before starting the proposed work, we recall some outstanding contributions by the authors in the classical sense; see \cite{benedetto2004wavelet,benedetto2004examples, albeverio2009multidimensional,kozyrev2002wavelet, jiang2004multiresolution}. The wavelet theory on local fields and analogous groups is discussed in \cite{benedetto2004wavelet,benedetto2004examples}. Further, the authors in \cite{albeverio2009multidimensional, kozyrev2002wavelet, jiang2004multiresolution} also explored the MRA technique and wavelets on the $p$-adic field $\mathbb{Q}_p$ of characteristic $p=0$. Furthermore, the authors \cite {jiang2004multiresolution} established the MRA technique on a local field in wavelet analysis along with interesting properties of wavelets. An example is also presented. From the work of Bahera and Jahan \cite{behera2012wavelet}, one can easily observe that to construct orthonormal wavelet packets in the space $L^{2}(K_q)$, {\it it only requires a single scaling function $\varphi$}. It is pretty similar to the construction of the orthonormal wavelet packets in $L^{2}({\mathbb R})$; see \cite{van2011wavelet} for a background. Likewise, we wish to construct orthonormal wavelet packets in $H^{\mathfrak s}(K_q)$. However, constructing wavelet packets is not a simple task in the Sobolev space, as a single scaling function is incapable of generating orthonormal functions at each level (due to the variant property of the Sobolev norm on dilation). This phenomenon gives a more general construction as scaling functions depend on the level. \\
\par {\bf Organization of the work}: In section~2, definitions and preliminaries of the local field $K_q$ are provided. In section~3, we have defined the space $H^{\mathfrak s}(K_q)$ and illustrated this space via seven distinct, exciting examples. Further, we stated some essential results and established the concept of the MRA technique in $H^{\mathfrak s}(K_q)$. Following this technique, we have constructed wavelet packets in the same space. Section 4 proves some important properties of the basic wavelet packets at the base level. In section 5, we have constructed wavelet packet functions at different scales generated by MRA in $H^{\mathfrak s}(K_q)$, and their orthogonality properties are derived. Further, an alternate technique is used to derive orthogonality properties of wavelet packet functions at different scales and construct three important examples (Haar wavelet packets, Weierstrass type wavelet packets, and 3-adic Cantor type wavelet packets). Finally, a conclusion is made in the last section. 
	
	
	\section{Definitions and preliminaries on local field}
	This section provides the definition of a local field $K_q$ and its important subsets, such as integer-ring, ideal, prime ideal, maximal ideal, group of units, fractional ideal, and their Haar measures, along with algebraic as well as topological proprieties.
	\begin{definition}[Non-Archimedean valued norm]
		If a mapping $|.|: K_q\rightarrow [0, \infty)$, where $q=p^{c}, p \ \text{a \ prime \ number},  c\in{\mathbb{N}}$ defined by $|a|$ satisfies the following axioms
		\begin{itemize}
			\item[(1)] For all $a\in K_q$, $|a| = 0$ if and only if $a=0$,
			\item[(2)] for all $a$ and $b$ belongs to $K_q$, $|ab|=|a|\ |b|$, 
			\item[(3)](Ultrametric inequality) For any two members $a$ and $b$ belongs to $K_q$ such that $|a|\neq|b| \implies |a+b|=\max\{|a|, |b| \}$.
			Then $|a|$ is called a non-Archimedean norm of $a$ in $K_q$, and the field $K_q$ is said to be a non-Archimedean valued field.  
		\end{itemize}
	\end{definition} 
	\begin{definition}[Local Field]
		\par A complete topological field which possesses the property of locally compact, totally disconnected, and non-Archimedean valued so-called \emph{a local field} $K_q$, for a prime $p$ and $c\in{\mathbb{N}}$ such that $q=p^{c}$. Note that $K_q^+$ and $K_q^{\ast}$ are two locally compact additive abelian and multiplicative groups in $K_q$.			
		\end{definition}
	{\bf Note:} The relevance of the integer ``q" lies in the structure and properties of the local field $K_q$. The integer ``q" in the context of local fields and non-Archimedean valued norms refers to the cardinality of the residue field, which has significant implications for the algebraic structure and properties of the local field. The residue field of a local field is obtained by considering the equivalence classes of elements in the local field modulo a maximal ideal. The cardinality (the number of elements) of the residue field is denoted by $q=p^{c}, p \ \text{a \ prime \ number},  c\in{\mathbb{N}}$.
		Now we are only interested in $K_q^+ (\equiv K_q )$ to define Haar measure and Haar integrals as follows:
	\subsection{Haar measure and integral on local field}
The Haar measure for a measurable set $E\subset K_q$ is defined by $\mu(E)\equiv |E|$ which satisfies invariance of translation (i.e., $|aE|=|Ea|=|E|$ for all $a \in K_q$). The quantity $d(\alpha a)=|\alpha|da$, where $|\alpha|$ is the \emph{valuation} of $\alpha$ for all $\alpha\neq 0 \in K_q$. 

The Haar integral for a Haar measurable function $\phi:K_q\rightarrow \mathbb{C}$ is defined by $\int_{K_q} \phi(x)dx$ satisfies translation invariance for all  $a \in K_q,\quad \int_{K_q} \phi(x+a)da=\int_{K_q} \phi(a)da$. Further, $|E|=\int_{K_q}\chi_E(a)da$, where $da$ is the Haar measure of $K_q$ and $\chi_E$ is the characteristic function of $E$.
	\subsection{Important subsets in local field}\label{sec2.2}
	Let us define the following sets on the local field as follows:
	\begin{itemize}
		\item[(1)] The set ${\mathcal D}_q=\{a\in K_q: |a|\leq 1\}$ denotes the \emph{ring of integers} in $K_q$, which is a unique maximal sub-ring, open, closed, as well as compact. The Haar measure of this set is $|{D_q}| = 1$. \\
		\item[(2)] The set ${\mathcal P}_q=\{a\in K_q:|a|<1\}$ denotes the \emph{prime ideal} in $K_q$. This ideal is a unique and maximal principle in $\mathcal{D}_q$. Moreover, it is an open, closed, and compact set, and its Haar measure is $|\mathcal{P}_q| = q^{-1}$. \\
		Before defining further important sets in $K_q$, let us see some useful facts in $K_q$. Since $K_q$ is totally disconnected. Therefore, the set $S=\{s^k: k\in\Z\}\cup \{0\}$ is a discrete set for some $s>0$ where the values of $|a|$ as $a$ runs over $K_q$ in discrete manner. Following this fact, there exists a fixed element $\mathfrak{p}$ of $\mathcal{P}_q$ with maximal absolute value, such that $\mathfrak{p}$ becomes a \emph{prime element} of $K_q$. Note that an ideal in $\mathcal{D}_q$ can be represented by $\mathcal{P}_q=\left\langle \mathfrak{p}\right\rangle = \mathfrak{p}\mathcal{D}_q$ with Haar measure $|\mathfrak{p}|=q^{-1}$.	
		\item[(3)] The set $\mathcal{D}_q^*=\mathcal{D}_q\setminus\mathcal{P}_q=\{a\in K_q: |a|=1\}$ denotes the \emph{group of units} in $K_q^*$. This set is open, closed, and compact. The Haar measure is $|\mathcal{D}_q^*| = 1-q^{-1}$.\\
		\item[(4)]  The set $\mathcal{P}_{q}^k=\{a\in K_q: |a|\leq q^{-k}\}=\mathfrak{p}^k\mathcal{D}_q, \ k\in\Z$ denotes the \emph{fractional ideals} in $K_q$. Each $\mathcal{P}_{q}^k$ is compact, open and is a subgroup of $K_q^+$ (see \cite{taibleson2015fourier}). The Haar measures are $|\mathcal{P}_{q}^k|=\mathfrak{p}^{k}=q^{-k}$, $k\in \Z$. If an element $a\neq 0$ belongs to $K_q$, then we can write $a=\mathfrak{p}^k a'$, where $a'\in\mathcal{D}_q^*$. 
	\end{itemize}	
	  Now from \cite{taibleson2015fourier}, we can define the quotient set in $K_q$, by $\mathcal{D}_q/\mathcal{P}_q \cong GF(q)$. Assume $\mathcal{U}=\{u_i\}_{i=0}^{q-1}$ is any fixed full set of coset representatives of $\mathcal{P}_q$ in $\mathcal{D}_q$. Now, if $a$ belongs $\mathcal{P}_{q}^k, k\in{\mathbb{Z}}$, then $a$ can be written in a unique manner by $a=\sum_{l=k}^{\infty} c_{l} \mathfrak{p}^{l}$, where $c_{l} \in \mathcal{U}$. 
	 
	 Let us generalize the Fourier transform on a local field $K_q$ with the help of Pontryagin duality. The additive group of $K_q^+$ of a local field $K_{q}$ is isomorphic to its character group $\Gamma$, and also $K_q^+$ is a self-dual (due to Pontryagin duality theorem). Following this fact, we can find a non-trivial, unitary, and continuous character $\chi \in \Gamma$ in such a way that it is always trivial on $\mathcal{D}_q$ but is non-trivial on $\mathcal{P}_{q}^{-1}$. Thus, this character $\chi$ can be obtained by starting with any non-trivial and rescaling character. We can define and classify this character $\chi$ in the following way, for $a,b\in K_q$ such that $\chi_b(a)=\chi(ba)$.	

	\begin{definition}\label{def1}
		Suppose $\phi$ is a member of $L^1(K_q)$, then the Fourier transform of the function $\phi$ is denoted by $\hat \phi$ and defined by
		\[
		\hat \phi(\xi)= \int_{K_q} \phi(x)\overline{\chi_{\xi}(x)}~{\rm d}x=\int_{K_q} \phi(x)\overline{\chi(\xi x)}~{\rm d}x = \int_{K_q} \phi(x)\chi(-\xi x)~{\rm d}x.
		\]
	\end{definition}
	where $\overline{\chi(\xi x)}=\chi(-\xi x)$. The following results can be easily proved on $K_q$ in a similar way on the real line $\mathbb{R}$ by using the Fourier analysis technique.
	\begin{itemize}
		\item[(1)] The map $\phi\rightarrow \hat \phi$ is a linear and bounded non-increasing norm transformation form $L^1(K_q)$ into $L^{\infty}(K_q)$, that is,
		 ${ (\phi_1+\phi_2)}^\wedge\rightarrow {\hat{\phi}_1}+{\hat{\phi}_2}$ and $\|\hat \phi\|_{L^{\infty}(K_q)}\leq \|\phi\|_{L^{1}(K_q)}$.
		\item[(2)] If $\phi \in L^1(K_q)$, then $\hat{\phi}(\xi)$ is uniformly continuous function on $\xi \in \Gamma$.
		\item[(3)] $\phi \in L^1(K_q) \cap L^2(K_q)$, then $\|\hat {\phi} \|_2=\|\phi\|_2$.
	\end{itemize}
	
	We borrow the definition of Fourier transform for a function $\phi$ in $L^2(K_q)$ provided in \cite{taibleson2015fourier} and let $\Phi_{\mathcal{P}_{q}^k}$ be the characteristic function of $\mathcal{P}_{q}^k$, for $k\in{\mathbb Z}$, then we define the following definition.
	\begin{definition}[Fourier transform of $L^2$-functions]
	Let $\phi_k=\phi\Phi_{\mathcal{P}_{q}^k}$. Then $\phi_k$ is the so called $k$-cut function of $\phi$, where $k\in \mathbb{Z}$. It is obvious that $\phi_k$ belongs to $L^1(K_q) \cap L^2(K_q)$. Then 
		\begin{eqnarray*}
	\hat{\phi}_k(\xi)=\int_{K_q} \phi_k(x)\overline{\chi_{\xi}(x)}~{\rm d}x=\int_{\left|x\right|\leq q^k} \phi(x)\overline{\chi_{\xi}(x)}~{\rm d}x.
	\end{eqnarray*}
	Hence, the Fourier transform of a function $\phi \in L^2(K_q)$ is provided by 
		\[
		\hat{\phi}(\xi)=\lim\limits_{k\rightarrow +\infty}\hat \phi_k(\xi)
		=\lim\limits_{k\rightarrow +\infty}\int_{\left|x\right|\leq q^k} \phi(x)\overline{\chi_{\xi}(x)}~{\rm d}x,\ {\text{a.e. for all} } \ x \in K_q,
		\]
		where the limit has to be understood in $L^2$-sense. 
	\end{definition}
\begin{definition}[Multiplication formula]
	If $\phi_1, \phi_2 \in L^2(K_q)$, then  	
\begin{eqnarray*}
\int_{K_q}{\hat{\phi}}_1(x){\phi}_2(x) {\rm d}x=\int_{K_q}{\phi}_1(x){\hat{\phi}}_2(x) {\rm d}x.
\end{eqnarray*}
\end{definition} 
	\begin{theorem}[Theorem 2.3 in~\cite{taibleson2015fourier}]
		The $L^2$-Fourier transform on $L^2(K_q)$ is unitary linear operator, i.e., the Parseval formula 
		\begin{eqnarray*}
			\int_{\Gamma}{\hat{\phi}}_1(\xi){\overline{\hat{\phi}}}_2(\xi) {\rm d}\xi=\int_{K_q}{\phi}_1(x){\overline{\phi}_2}(x) {\rm d}x.
		\end{eqnarray*}
	\end{theorem}
	
We would choose a character $\chi_{\lambda}\in \mathcal{D}_q$ of $K_q$, such that it is trivial on $\mathcal{D}_q$, i.e., $\chi_{\lambda}|_{\mathcal{D}_q}\equiv 1$. Also, as character on $\mathcal{D}_q, \chi_{\lambda_1} = \chi_{\lambda_2}$ if and only if ${\lambda_1-\lambda_2}\in \mathcal{D}_q$. That is, $\chi_{\lambda_1}=\chi_{\lambda_2}$ if $\lambda_1+\mathcal{D}_q=\lambda_2+\mathcal{D}_q$ and $\chi_{\lambda_1}\neq \chi_{\lambda_2}$ if $(\lambda_1+\mathcal{D}_q)\cap (\lambda_2+\mathcal{D}_q)=\emptyset$. Hence, if $\{\lambda(n)\}_{n=0}^{\infty}$ is a complete list of distinct coset representative of $\mathcal{D}_q$ in $K_q^+$ implies that the list $\{\chi_{\lambda(n)}\}_{n=0}^{\infty}$ is distinct characters on $\mathcal{D}_q$. The completeness of this list is proved in~\cite{taibleson2015fourier} and we have the following proposition on the list $\{\chi_{\lambda(n)}\}_{n=0}^{\infty}$.
	\begin{proposition}\label{pro:1}
		Suppose $\{\lambda(n)\}_{n=0}^{\infty}$ is a complete list of distinct coset representatives of $\mathcal{D}_q$ in $K_q$. Then $\{\chi_{\lambda(n)}\}_{n=0}^{\infty}$ is a complete list of distinct characters on $\mathcal{D}_q$. Further, this list $\{\chi_{\lambda(n)}\}_{n=0}^{\infty}$ is a complete orthonormal system on $\mathcal{D}_q$.
	\end{proposition}
	By using this list of characters $\{\chi_{\lambda(n)}\}_{n=0}^{\infty}$, we can  define the Fourier coefficients of $\phi\in L^1(\mathcal{D}_q)$ as follows
	\[
	\hat{\phi}(\lambda(n))=\int_{\mathcal{D}_q}\phi(x)\overline{\chi_{\lambda(n)}(x)}{\rm d}x.
	\]
	The Fourier series $\sum\limits_{n=0}^{\infty}\hat{\phi}(\lambda(n))\chi_{\lambda(n)}(x)$ converges to $\phi$ in $L^2(\mathcal{D}_q)$ and following relation holds
	\[
	\int_{\mathcal{D}_q}|\phi(x)|^2{\rm d}x= \sum\limits_{n=0}^{\infty}|\hat{\phi}(\lambda(n))|^2.
	\]
	Further, $\forall n\in\N_0$, $\phi\in L^1(\mathcal{D}_q)$ and $\hat \phi(\lambda(n))=0$, then $\phi=0$ a.e..
	
The above results hold irrespective of the character's ordering. We now begin to impose a natural ordering on this sequence $\{\lambda(n)\}_{n=0}^{\infty}$. As $\mathcal{P}_q$ is the prime ideal in $\mathcal{D}_q$, $\mathcal{D}_q/\mathcal{P}_q\cong GF(q)$, and $\rho : \mathcal{D}_q \rightarrow \Gamma$ the canonical homomorphism of $\mathcal{D}_q$ on to $\Gamma$. As we know that $\Gamma=GF(q)$ is a $c$-dimensional vector space over the Galois field $GF(p)\subset \Gamma$. Now, we consider a set $\{1=e_0, e_1, e_2, \cdots, e_{c-1}\} \subset\mathcal{D}_q^*$ in such a way that $\{\rho(e_j)\}_{j=0}^{c-1}$ is a basis of $GF(q)$ over $GF(p)$. If $n\in \N_0$ such that $0\leq n< q$, then $n$ can be written by
	\[
	n=u_0+u_1 p+\cdots+u_{c-1} p^{c-1},\ \text{for}\ 0\leq u_k<p, \ \text{and}\
	 k=0,1,\cdots,c-1.
	\]
	Define
	\begin{equation}\label{e.undef}
	\lambda(n)=(u_0+u_1 e_1+\cdots+u_{c-1}e_{c-1})\mathfrak{p}^{-1}.
	\end{equation}
	Now, write
	\[
	n=v_0+v_1q+v_2q^2+\cdots+v_sq^{s},\ \text{for}\ 0\leq v_k<q,\ \text{and}\ k=0,1,2,\cdots,s,
	\]
	and define
	\[
	\lambda(n)=\lambda(v_0)+\mathfrak{p}^{-1}\lambda(v_1)+\cdots+\mathfrak{p}^{-s}\lambda(v_s).
	\]
	Note that $\lambda(0)=0$ and $\{\lambda(n)\}_{n=0}^{q-1}$ is a complete set of coset representatives of $\mathcal{D}_q$ in $\mathcal{P}_q^{-1}$ (see~\cite{taibleson2015fourier}). Hence, $\{\lambda(n)\mathfrak{p}\}_{n=0}^{q-1}$ is a complete set of coset representatives of $\mathcal{P}_q$ in $\mathcal{D}_q$. Thus
	\[
	\{\lambda(n){\mathfrak{p}}\}_{n=0}^{q-1}\cong \mathcal{D}_q/\mathcal{P}_q\cong GF(q)\cong \{\rho(e_j)\}_{j=0}^{c-1}.
	\]
	
Note that $\lambda(m+n)=\lambda(m)+\lambda(n)$ is not true in general. But
	\begin{equation}\label{eq.un}
	\lambda(rq^k+s)=\lambda(r)\mathfrak{p}^{-k}+\lambda(s)\quad{\rm if}~r\geq 0, k\geq 0~{\rm and}~0\leq s <q^k.
	\end{equation}
	We shall denote $\chi_n=\chi_{\lambda(n)}$ for all $n\geq 0$ and we know that the sequence $\{\chi_n\}_{n=0}^{\infty}$ is a complete set of characters on $\mathcal{D}_q$. Let $\mathcal{U}=\{u_i\}_{i=0}^{q-1}$ be a fixed set of coset representatives of $\mathcal{P}_q$ in $\mathcal{D}_q$. Then every $x\in K_q$ can be written in a unique way as follows
	\[
	x=x_0+\sum\limits_{k=1}^n v_k\mathfrak{p}^{-k}, \quad x_0\in \mathcal{D}_q, v_k\in \mathcal{U}.
	\]
	
	Now, we can represent a character $\chi$ on $K_q$ by
	\begin{equation}\label{chi}
	\chi(e_{\mu}\mathfrak{p}^{-j})=
	\left\{
	\begin{array}{lll}
	\exp(\frac{2\pi i}{p}), & \mu=0~\mbox{and}~j=1,\\
	1, & \mu=1,\cdots,c-1~\mbox{or}~j\neq 1.
	\end{array}
	\right.
	\end{equation}
	
	To construct MRA, basic wavelets, and wavelet packet functions at different scales in $H^{\mathfrak s}(K_q)$. We must take special care of dilation as the single scaling function does not work here. To do this, we require analogous notions of translation and dilation. Note that $\bigcup\limits_{j\in\Z}\mathfrak{p}^{-j}\mathcal{D}_q=K_q$, then $\mathfrak{p}^{-1}$ can be used as the dilation ($|\mathfrak{p}^{-1}|=q$) and since $\{\lambda(n): n\in\N_0\}$ is a complete list of distinct coset representatives of $\mathcal{D}_q$ in $K_q$, therefore, this set $\{\lambda(n): n\in\N_0\}$ can serve as the translation set. 
		
	\section{Multiresolution analysis in $H^{\mathfrak s}(K_q)$}
	In this section, we define the Sobolev space ($H^{\mathfrak s}(K_q)$, where ${\mathfrak s}\in {\R}$) and construct MRA and wavelet packets in this space. Examples, including Weierstrass type and Cantor type functions, illustrate the Sobolev space and its special cases. Further, some remarkable results, including orthogonal properties, are derived.
	\begin{definition}[The space of tempered distributions $\mathcal{S'}(K_q)$] A tempered distribution is a continuous linear map $\mathcal{S}(K_q)\mapsto \mathbb{C}$, where the space $\mathcal{S}(K_q)$ denotes Schwartz space. In other words: If $f$ is a locally integrable function of polynomial growth   on $K_q$ then $f$ it generates a tempered  distribution in $\mathcal{S'}(K_q)$ as follows: 
		\begin{eqnarray*}
		\bigl\langle f, \phi \bigr\rangle =  \int_{K_q} f(x) \phi(x) {\rm d}x, \ \forall \phi \in \mathcal{S}(K_q),
		\end{eqnarray*}
where $\mathcal{S}(K_q)$ and $\mathcal{S'}(K_q)$ denote Schwartz space and its dual, respectively. 	
			\end{definition}
	
	\begin{definition}[Sobolev space $H^{\mathfrak s}(K_q)$]\label{def2}
		It consists of all tempered distributions $\phi$ in $\mathcal{S'}(K_q)$ such that
		\begin{equation*}
		\left\|\phi\right\|_{H^{\mathfrak s}(K_q)}=\int_{K_q}(1+|\xi|^2)^{\mathfrak s} |\hat{\phi}(\xi)|^2 {\rm d}\xi< +\infty.
		\end{equation*}
		Using the Parseval formula and multiplication formula, the corresponding inner product can be written as follows 
		\begin{eqnarray*}
			\bigl\langle \phi_1, \phi_2 \bigr\rangle_{H^{\mathfrak s}(K_q)} =  \int_{K_q} (1+|\xi|^2)^{\mathfrak s} \hat{\phi_1}(\xi)\overline{\hat{\phi_2}(\xi)}{\rm d}\xi,
		\end{eqnarray*}
		where $\hat{\phi}$ is defined in Definition \ref{def1}.
	\end{definition}

We want to illustrate the Sobolev space $H^{\mathfrak s}(K_q)$ and its subspace via examples, including fractals. \begin{example}\label{ex1}
	Let $\theta >-1$, and $x\in K_q$, then the function $\phi(x) $ defined by 
	\begin{eqnarray}\label{eqe1}
		\phi(x) =
		\left\{
		\begin{array}{ll}
			{|x|}^{\theta}  & \mbox{if } |x| \leq 1 \\
			0 & \mbox{if } |x| >1,
		\end{array}
		\right.
	\end{eqnarray}
	belongs to $H^{\mathfrak s}(K_q)$ for every real value of $\mathfrak s <\theta+1/2$. 
\end{example}
\begin{solution} The Fourier transform of the given function $\phi(x)$ is as follows:
	
	\begin{eqnarray}
		\hat{\phi}(\xi) =
		\left\{
		\begin{array}{ll}
			\frac{1-{q}^{-1}}{1-{q}^{-(1+\theta)}}  & \mbox{if } |\xi| \leq 1 \\
			\frac{1-{q}^{-1}}{1-{q}^{-(1+\theta)}} \frac{1}{{|\xi|}^{1+\theta}} & \mbox{if } |\xi| >1,
		\end{array}
		\right.
	\end{eqnarray}
	then 
	\begin{eqnarray*}
		\left\|\phi\right\|_{H^{\mathfrak s}(K_q)}
		&\leq& \int_{|\xi|\leq 1}(1+|\xi|^2)^{\mathfrak s} \frac{(1-{q}^{-1})^2}{(1-{q}^{-(1+\theta)})^2} {\rm d}\xi+2^{\mathfrak s}\int_{|\xi|>1}|\xi|^{2\mathfrak s-2\theta-2} \frac{(1-{q}^{-1})^2}{(1-{q}^{-(1+\theta)})^2}	 {\rm d}\xi<\infty,
	\end{eqnarray*}
	since the first integral on the right-hand side of the above expression is bounded and the second integral is convergent for every real value of $\mathfrak{s}<\theta+1/2$.
\end{solution}
\begin{example}\label{ex2}
	For $x\in K_q$, then the function $\phi(x) $ defined by 
	\begin{eqnarray}\label{eqe2}
		\phi(x) =
		\left\{
		\begin{array}{ll}
			{ln\frac{1}{|x|}}  & \mbox{if } |x| \leq 1\\
			0 & \mbox{if } |x| >1,
		\end{array}
		\right.
	\end{eqnarray}
	belongs to $H^{\mathfrak s}(K_q)$ for every real value of $\mathfrak s <1/2$. 
\end{example}
\begin{solution} The Fourier transform of the given function $\phi(x)$ is as follows:
	
	\begin{eqnarray}
		\hat{\phi}(\xi) =
		\left\{
		\begin{array}{ll}
			\frac{\ln{q}}{1-{q}^{-1}}{q}^{-1}  & \mbox{if } |\xi| \leq 1 \\
			\frac{\ln{q}}{1-{q}^{-1}}{q}^{-1} \frac{1}{{|\xi|}} & \mbox{if } |\xi| >1,
		\end{array}
		\right.
	\end{eqnarray}
	then 
	\begin{eqnarray*}
		\left\|\phi\right\|_{H^{\mathfrak s}(K_q)}
		&\leq& \int_{|\xi|\leq 1}(1+|\xi|^2)^{\mathfrak s} \frac{(\ln{q})^2}{(1-{q}^{-1})^2}{q}^{-2} {\rm d}\xi+2^{\mathfrak s}\int_{|\xi|>1}|\xi|^{2\mathfrak s-2} \frac{(\ln{q})^2}{(1-{q}^{-1})^2}{q}^{-2} {\rm d}\xi<\infty,
	\end{eqnarray*}
since the first integral on the right-hand side of the above expression is bounded and the second integral is convergent for every real value of $\mathfrak{s}<1/2$.
\end{solution}
\begin{example}\label{ex3}
	Let us recall the definition of the fractional ideas defined in section \ref{sec2.2} as $\mathcal{P}_{q}^k=\{a\in K_q: |a|\leq q^{-k}\}=\mathfrak{p}^k\mathcal{D}_q, \ k\in\Z$ in $K_q$ and characteristic function of $\mathcal{P}_{q}^k$:
	
	\begin{eqnarray}\label{eqe3}
		\Phi_{\mathcal{P}_{q}^k}
		(x) =
		\left\{
		\begin{array}{ll}
			{1 } & \mbox{if } x \in {\mathcal{P}_{q}^k} \\
			0 & \mbox{if } x \in K_q \setminus {\mathcal{P}_{q}^k},
		\end{array}
		\right.
	\end{eqnarray}
	belongs to $H^{\mathfrak s}(K_q)$ for every real value of $\mathfrak s$. 
\end{example}
\begin{solution} The Fourier transform of the given characteristic function $\Phi_{\mathcal{P}_{q}^k}
	(x)$ is as follows:
	
	\begin{eqnarray}
		\hat{\Phi}_{\mathcal{P}_{q}^k}(\xi) =
		\left\{
		\begin{array}{ll}
			{q}^{-k}  & \mbox{if } |\xi| \leq {q}^{k}\\
			0 & \mbox{if } |\xi| > {q}^{k},
		\end{array}
		\right.
	\end{eqnarray}
	then 
	\begin{eqnarray*}
		\left\|\Phi_{\mathcal{P}_{q}^k}
		\right\|_{H^{\mathfrak s}(K_q)}
		&=&\int_{|\xi|\leq {q}^{k}}(1+|\xi|^2)^{\mathfrak s} {q}^{-2k} {\rm d}\xi<\infty,	
	\end{eqnarray*}
	since the integral on the right-hand side of the above expression is bounded for every real value of ${\mathfrak{s}}$.
\end{solution}
\begin{example}\label{ex4}
	Let $0<\theta, \vartheta, \theta+\vartheta <1$, and $x\in K_q$, then the function $\phi(x) $ defined by 
	\begin{eqnarray}\label{eqe4}
		\phi(x) =
		\left\{
		\begin{array}{ll}
			\frac{1-q^{-(\theta+\vartheta)}}{1-q^{\theta+\vartheta-1}}{|x|}^{\theta+\vartheta-1}  & \mbox{if } |x| \leq 1 \\
			0 & \mbox{if } |x| >1,
		\end{array}
		\right.
	\end{eqnarray}
	belongs to $H^{\mathfrak s}(K_q)$ for every real value of $\mathfrak{s}<\theta+\vartheta-1/2$. 
\end{example}
\begin{solution} The Fourier transform of the given function $\phi(x)$ is as follows:
	
	\begin{eqnarray}
		\hat{\phi}(\xi) =
		\left\{
		\begin{array}{ll}
		 \frac{1-q^{-1}}{1-q^{-(\theta+\vartheta)}} & \mbox{if } |\xi| \leq 1 \\
		\frac{1-q^{-1}}{1-q^{-(\theta+\vartheta)}}		|\xi|^{-(\theta+\vartheta)} & \mbox{if } |\xi| >1.
		\end{array}
		\right.
	\end{eqnarray}
	Then 
	\begin{eqnarray*}
		\left\|\phi\right\|_{H^{\mathfrak s}(K_q)}
		&\leq& \int_{|\xi|\leq 1}(1+|\xi|^2)^{\mathfrak s} \frac{(1-{q}^{-1})^2}{(1-{q}^{-(\theta+\vartheta)})^2} {\rm d}\xi+2^{\mathfrak s}\int_{|\xi|>1}|\xi|^{2(\mathfrak s-\theta-\vartheta)} \frac{(1-{q}^{-1})^2}{(1-{q}^{-(\theta+\vartheta)})^2}	 {\rm d}\xi<\infty,
	\end{eqnarray*}
since the first integral on the right-hand side of the above expression is bounded and the second integral is convergent for every real value of $\mathfrak{s}<\theta+\vartheta-1/2$.
\end{solution}
\subsection{More fractal-like examples on a special case of Sobolev space}
This subsection illustrates special cases of Sobolev space on a local field of finite characteristic $p>0$ with the help of examples having fractal-like structures. Here, we will be showing three more examples that belong to special cases of Sobolev spaces. 

\begin{example}\label{ex5}[Weierstrass type function]
The well-known Weierstrass function on the real line is a typical fractal. In the same connection, we are interested in illustrating the special case of Sobolev space by the example of the Weierstrass type function.

	Let $p=2$ and $c=1$ such that $q=p^c=2$, then $K_q$ reduces to $K_2$. Hence, the subset of the Sobolev space is defined by $H^{\mathfrak{s}}(K_2)$. For $x\in K_2$, $x=\sum_{j=s}^{+\infty} {2}^{-j} x_{j} $, $x=\{0,1\}, j=s,s+1,s+2,\dots,-2,-2,0,1,2,\dots$, with $|\mathfrak{p}|=1/2$. then the Weierstrass type function ${\mathcal W}(x) $ defined by 
	\begin{eqnarray}\label{eqe5}
		{\mathcal W}_2(x) =
		\left\{
		\begin{array}{ll}
			\sum_{j=1}^{+\infty} {2}^{-j} x_{j}   & \mbox{if } x\in {\mathcal P}^{1}_{2} \\
			0 & \mbox{if }  x\notin {\mathcal P}^{1}_{2},
		\end{array}
		\right.
	\end{eqnarray}
	belongs to $H^{\mathfrak s}(K_2)$ for every real value of $\mathfrak s <-1/2$. 
\end{example}
\begin{solution} The Fourier transform of the Weierstrass type function ${\mathcal W}_2(x)$ is as follows:
	
	\begin{eqnarray}\label{eqe5.1}
		\hat{{\mathcal W}_2}(\xi) = 1/4, \ \ \ \mbox{where}\ \ |\xi|=2^{l},\ \  l\in {\mathbb{Z}} \ \ \mbox{such that}\ \   l\leq 1. 
		\end{eqnarray}
	Then 
	\begin{eqnarray*}
		\left\|{\mathcal W}_2\right\|_{H^{\mathfrak s}(K_q)}
		&=&\frac{1}{16}\int_{K_2}(1+|\xi|^2)^{\mathfrak s} {\rm d}\xi <\infty,
	\end{eqnarray*}
	since the last integral on the right-hand side of the above expression is convergent for every real value of $\mathfrak{s}<-1/2$.
\end{solution}
\begin{example}\label{ex6}[3-adic Cantor type function (ladder of devil)]
		The well-known 3-adic Cantor function on the real line is a typical fractal. In the same connection, we are interested in illustrating the special case of Sobolev space by the example of the 3-adic Cantor type function.
	
	Let $p=3$ and $c=1$ such that $q=p^c=3$, then $K_q$ reduces to $K_3$. Hence, the subset of the Sobolev space is defined by $H^{\mathfrak{s}}(K_3)$. For all $x\in K_3$, $x=\sum_{j=-s}^{+\infty} {\mathfrak p}^{j} x_{j}$, $x=\{0,1,2\}, j=-s,-s+1,-s+2,\dots,-2,-2,0,1,2,\dots$, with $|\mathfrak{p}|=1/3$ such that support of ${\mathcal C}(x)={\mathcal P}^1_3$ and $\forall x\in{\mathcal P}^1_3$ with $x=\sum_{j=0}^{+\infty} {3}^{-j} x_{j} $, $x=\{0,1,2\}, j=0,1,2,\dots$, then the 3-adic Cantor type function ${\mathcal C}(x)$ defined by 
	\begin{eqnarray*}\label{eqe6}
		{\mathcal C}(x) =
		\left\{
		\begin{array}{ll}
			\sum_{j=0}^{k-2} {2}^{-j-1}  (x_{j}-1) +{2}^{-k} & \mbox{if }  \ \exists\ k\geq 1; \ni x_j\neq 0 \ ( 0\leq j\leq k-2); \ x_{k-1}=0\\			
			\sum_{j=0}^{k-2} {2}^{-j-1} (x_{j}-1)  +{2}^{-k}& \mbox{if}\ x_j\neq \ 0,\  \forall \ 0\leq j<+\infty\\
			0 & \mbox{if }  x\notin {\mathcal P}^{1}_{3},
		\end{array}
		\right.
	\end{eqnarray*}
	belongs to $H^{\mathfrak s}(K_3)$ for every real value of $\mathfrak s <-1/2$. 
\end{example}
\begin{solution} The Fourier transform of the 3-adic Cantor type function ${\mathcal C}(x)$ is as follows:
	
	\begin{eqnarray*}
		\hat{{\mathcal C}}(\xi) = 1/2, \ \ \ \mbox{where}\ \ |\xi|=3^{l},\ \  l\in {\mathbb{Z}} \ \ \mbox{such that}\ \   l\leq 0. 
	\end{eqnarray*}
	Then 
	\begin{eqnarray*}
		\left\|{\mathcal C}\right\|_{H^{\mathfrak s}(K_3)}
		&=&\frac{1}{4}\int_{K_3}(1+|\xi|^2)^{\mathfrak s} {\rm d}\xi<\infty,
	\end{eqnarray*}
	since the integral on the right-hand side of the above expression is convergent for every real value of $\mathfrak{s}<-1/2$.
\end{solution}
\begin{example}\label{ex7}
	As $q=p^c$, if we take $c=1$, then $K_q$ reduces to $K_p$. Let $\theta \in \mathbb{R}$, and $x\in K_p$, then the function $\kappa_{\theta}(x) $ defined by 
	\begin{eqnarray}\label{eqe7}
		\kappa_{\theta}(x) =
		\left\{
		\begin{array}{ll}
			\left\{\left(\frac{1-p^{\theta}}{1-p^{\theta-1}}\right)\pi_{-\theta(x)} +\left(\frac{1-p^{\theta}}{1-p^{\theta+1}}\right)\right\}\Phi_{{\mathcal P}^{0}_{p}(x)} & \mbox{if } \theta \neq 0,-1 \\
			\delta(x)  & \ \mbox{if}\  \theta = 0\\
			\left(1-1/p\right)(1-\log_p|x|)\Phi_{{\mathcal P}^{0}_{p}}(x) & \mbox{if } \ \theta= -1,
		\end{array}
		\right.
	\end{eqnarray}
	belongs to $H^{\mathfrak s}(K_q)$ for every real value of $\mathfrak s <-(1+2\theta)/2$. Where the distribution $\pi_{-\theta}$ is defined as follows
	\begin{eqnarray*}
		<\pi_{\theta},\varphi>=\int_{K_p} \frac{(\varphi(x)-\varphi(0))}{|x|^{1-\theta}}{\rm d}x, \ \ \forall \varphi \in {\mathcal S}(K_p) \ \mbox{and} \ \theta (\neq 0) \in {\mathbb{R}}.
	\end{eqnarray*}
\end{example}
\begin{solution} The Fourier transform of the given function $\kappa_{\theta}(x)$ is as follows:
	
	\begin{eqnarray}\label{eqe7.1}
		\hat{\kappa_{\theta}}(\xi) =(1+|\xi|^2)^{\theta}	\end{eqnarray}
	Then 
	\begin{eqnarray*}
		\left\|\phi\right\|_{H^{\mathfrak s}(K_2)}
		&=&\int_{K_p}(1+|\xi|^2)^{\mathfrak s} (1+|\xi|^2)^{2\theta} {\rm d}\xi<\infty,
		\end{eqnarray*}
	since the integral on the right-hand side of the above expression is convergent for every real value of $\mathfrak{s}<-2\theta-1/2$.
\end{solution}

	\subsection{Multiresolution analysis in $ H^{\mathfrak s}(K_q)$}
	
	The theory of MRA, including wavelet packets in $L^2(K_q)$, was developed by Behera and Jahan \cite{behera2012wavelet}. In this subsection, we construct MRA in $H^{\mathfrak s}(K_q)$ in the following way:	
	\begin{definition}
		 The MRA of $H^{\mathfrak s}(K_q)$ is a sequence $\{\mathsf{V}_j\}_{j\in\mathbb{Z}}$ of closed linear subspaces of $H^{\mathfrak s}(K_q)$ such that 
		\begin{enumerate}
			\item[(1)] $\mathsf{V}_j\subset \mathsf{V}_{j+1}$; 
			\item[(2)] $\overline{\bigcup\limits_{j=-\infty}^{+\infty}\mathsf{V}_j}=H^{\mathfrak s}(K_q)$ and $\bigcap_{j=-\infty}^{+\infty}\mathsf{V}_j= \{0\}$;
			\item[(3)] $\phi \in \mathsf{V}_j$ if and only if $\phi(\mathfrak{p}^{-1}\cdot)\in \mathsf{V}_{j+1}$ for all $j\in\Z$;
			\item[(4)]for every $ j$, there is a function $\varphi^{(j)}\in \mathsf{V}_0$, called the \emph{scaling function}, such that $\{\varphi^{(j)}(\cdot-\lambda(k)): k\in\N_0\}$ forms an orthonormal basis for $\mathsf{V}_0$.
		\end{enumerate}
	\end{definition}
	\begin{proposition} \label{sec:pro:1}
		If ${\mathfrak s}\in {\mathbb{R}}$, $j\in \mathbb{Z}$ and $\varphi^{(j)}\in H^{\mathfrak s}(K_q)$, then functions $\varphi^{(j)}_{j,k}(\cdot)=q^{j/2} \varphi^{(j)}(\mathfrak{p}^{-j}\cdot-\lambda(k)),
		\ k\in {\N_{0}}$ are orthonormal in $H^{\mathfrak s}(K_q)$ if and only if
		\begin{eqnarray}
		\sum_{k=0}^{\infty} (1+p^{-2j} (|\xi|+\lambda(k))^2)^{\mathfrak s} |\hat{\varphi}^{(j)}(\xi+\lambda(k))|^2=1 \label{eq:1.1}
		\end{eqnarray}
		for a.e. $\xi \in K_q$. Following this, we get the bound below
		\begin{eqnarray*}
			|\hat{\varphi}^{(j)}(\mathfrak{p}^{j}\xi)|\leq (1+|\xi|^2)^{-{\mathfrak s}/2}.
		\end{eqnarray*}
	\end{proposition}
	\begin{proof}
		Proof of this proposition is easy; we omit it for the reader.
	\end{proof}
	For brevity, we shall write $\varphi^{(j)}_{j,k}$ instead of $\varphi_{j,k}$, which indicates that we are handling the $j^{th}$ level wavelets.
	The notation  $\mathsf{W}_j$ denotes the orthogonal complement of $\mathsf{V}_j$ in $\mathsf{V}_{j+1}$.
	\begin{proposition}\label{sec:pro:2}
		For every
		$j\in {\mathbb{Z}}$, let $\{\varphi^{(j)}\}$ be a sequence of elements of $H^{\mathfrak s}(K_q)$ such that the distributions $\varphi^{(j)}_{j,k}(\cdot)=q^{j/2} \varphi^{(j)}({\mathfrak{p}^{-j} \cdot-\lambda(k)}), k\in {\N_{0}}$, are orthonormal
		in $H^{\mathfrak s}(K_q)$. If $P_j$ is the orthogonal projection from $H^{\mathfrak s}(K_q)$ onto
		$\mathsf{V}_{j}=\overline{>\varphi^{(j)}_{j,k}: k\in \N_{0}<}$, then for each $h\in H^{\mathfrak s}(K_q)$, we have
		\begin{eqnarray*}
			\lim_{j\rightarrow + \infty} \left( \left\| P_{j} h\right\|_{H^{\mathfrak s}(K_q)}
			\int_{K_q} (1+|\xi|^2)^{2s} |\hat{h}(\xi)|^2 |\hat{\varphi}^{(j)} (\mathfrak{p}^{j} \xi)|^2 {\rm d}\xi\right)=0.
		\end{eqnarray*}
		Further, if there are $A$, and $\theta>0$, such that
		\begin{eqnarray*}
			\int_{K_q} (1+|\xi|^2)^{\theta} |\hat{\varphi}^{(j)} (\xi)|^2 {\rm d}\xi \leq A
		\end{eqnarray*}
		for each $j\leq 0$, then $\bigcap_{j=-\infty}^{+\infty} \mathsf{V}_j=\{0\}$.
	\end{proposition}
	\begin{proof}
		It is also not hard to prove this proposition; we omit it.
	\end{proof}
	\subsection{Construction of wavelets in $H^{\mathfrak s}({K_q})$}\label{sec:3.2}
		The essential results of the previous subsection will be used to construct wavelets in $H^{\mathfrak s}({K_q})$. 
\begin{definition}[Integral-periodic]
	We consider a function $\phi$ on $K_q$ \emph{integral-periodic} if
	\[
	\phi(x+\lambda(k))=\phi(x)~\mbox{for all}~k\in\N_0.
	\]
	\begin{itemize}
		\item[(1)] $\chi_k(\lambda(l))=\chi(\lambda(k)\lambda(l))=1$ for all $k, l\in\N_0$.
		\item[(2)] The function $m_0$ is integral-periodic.
	\end{itemize}
	The above facts (1) and (2) were proved in~\cite{jiang2004multiresolution}.
\end{definition}

	\begin{enumerate}
		\item Following the additional density condition below 
		\begin{eqnarray*}
			\lim_{j\rightarrow +\infty} |\hat{\varphi}^{(j)} (\mathfrak{p}^{j} \xi)|=(1+|\xi|^2)^{-{\mathfrak s}/2},
		\end{eqnarray*}
		 and from Proposition \ref{sec:pro:2}, since
		\begin{eqnarray*}
			\left\|g-P_{j}g\right\|_{H^{\mathfrak s}(K_q)}=\left\|g\right\|_{H^{\mathfrak s}(K_q)} - \left\|P_{j}g\right\|_{H^{\mathfrak s}(K_q)} \rightarrow 0,
		\end{eqnarray*}
		as $j\rightarrow + \infty$, for every $g\in H^{\mathfrak s}(K_q)$. Therefore, $\overline{\bigcup_{j=-\infty}^{+\infty}{\mathsf{V}_{j}}}$'s is dense in $H^{\mathfrak s}(K_q)$.
		
		\item By definition of MRA, $\mathsf{V}_{j}$ is the set of all $\phi\in H^{\mathfrak s}(K_q)$ such that
		\begin{eqnarray*}
			\hat{\phi}(\xi)=m(\mathfrak{p}^{j}\xi) \hat{\varphi}^{(j)}(\mathfrak{p}^{j}\xi),
		\end{eqnarray*}
		where $m\in L^{2}(\mathcal{D}_q)$ is integral-periodic.
	 Following this fact, we can find the Fourier transform of $q^{j/2} \varphi^{(j)}(\mathfrak{p}^{-j}\cdot-\lambda(k))$ is
		$ q^{-j/2}\overline{\chi_k(\mathfrak{p}^j\xi)} {\hat{\varphi}}^{(j)}(\mathfrak{p}^{j} \xi)$.
		
		\item We have $\mathsf{V}_{j} \subset \mathsf{V}_{j+1}$ for every $j\in {\mathbb{Z}}$ if and only if there exist periodic integral functions
		$m_{0}^{(j)}\in L^{2}({\mathcal{D}_q})$ such that the following result holds:
		\begin{eqnarray*}
			\hat{\varphi}^{(j)}(\xi)=m_{0}^{(j+1)}(\mathfrak{p}\xi) \hat{\varphi}^{(j+1)}(\mathfrak{p}\xi).
		\end{eqnarray*}
		
	\end{enumerate}
	Further, we know that $\varphi^{(j)}$ and $\varphi^{(j+1)}$ satisfy the hypothesis of Proposition \ref{sec:pro:1}, and these filters
	satisfy condition \ref{filter0} (a) given in section \ref{sec4}.
	
	\section{Basic wavelet packets in $H^{\mathfrak s}(K_q)$ space at base level}\label{sec4} 
	In this section, we focus mainly on deriving important results of basic wavelet packets at the base level in $H^{\mathfrak s}(K_q)$. Let $\{{\rm e}_t: t\in {\N_{0}}\}$ be a basis in $H^{\mathfrak s}({K_q})$. We can divide this system into $q$-parts by assuming $\{{\rm e}_{t,l}:t\in \N_{0},0\leq l\leq q-1\}$ in $\ell^2({\N_{0}})$. In general consider $q$-sequences (of $j^{th}$ level wavelet) $\{\alpha^{(j)}_{k,l}:k\in \N_{0},0\leq l\leq q-1, j\in \Z\}$. Define
	the elements $\{\phi^{(j)}_{k,l}: k\in \N_{0},0\leq l\leq q-1, j\in {\Z}\} \in H^{\mathfrak s}(K_q)$ by
	\begin{eqnarray}
	\begin{aligned}
	\phi^{(j)}_{k,l}=q^{1/2} \sum_{t\in {\N_{0}}} \alpha^{(j)}_{2k-t,l} {\rm e}_{t,l}; \forall \ l=0,1,2,...,q-1. \end{aligned} \label{eq:2.1}
	\end{eqnarray}
	
	We will characterize the systems $\{\alpha^{(j)}_{k,l}:k\in \N_{0},0\leq l\leq q-1, j\in \Z\}$ for which $\{\phi^{(j)}_{k,l}: k\in \N_{0},0\leq l\leq q-1, j\in {\Z}\}$ is still an orthonormal system in $H^{\mathfrak s}(K_q)$. To
	accomplish this, we assume
	
	\begin{eqnarray}
	\begin{aligned}
	&&m_l^{(j)}(\xi)= \sum_{k\in \N_{0}} \alpha^{(j)}_{k,l} \hat{e}_{t,l}.
	\end{aligned} \label{eq:2.2}
	\end{eqnarray}
	
	These functions belongs to $L^2(\mathcal{D}_q)$ and they are integral-periodic, since the sequences $\{\alpha^{(j)}_{k,l}:k\in \N_{0},0\leq l\leq q-1, j\in \Z\}$ are in $l^2({\N_{0}})$.
	Now, we consider the matrix
	
	\begin{eqnarray}
	M^{(j)}(\xi)=m^{(j)}_{l} \left(\mathfrak{p}\xi+\mathfrak{p}\lambda(k)\right)^{q-1}_{l,k=0}
	, \quad \xi \in {K_q}. \label{eq:2.3}
	\end{eqnarray}
\begin{proposition} \label{sec:pro:3}
		The matrix $M^{(j)}(\xi)$ given by \eqref{eq:2.2} and \eqref{eq:2.3} is unitary if and only if
		\begin{eqnarray*}
			&& (i) \ \sum_{t\in {\N_{0}}} \alpha^{(j)}_{t-2k,l} \overline{\alpha^{(j)}_{t,l}}=\frac{1}{q} \delta_{k,0}^{l}, \ \forall \ l=0,1,2,...,q-1.\\
			&& (ii) \ \sum_{t\in {\N_{0}}}\alpha^{(j)}_{t-2k,l} \overline{\alpha^{(j)}_{t,r}}=0, \ \forall \  0\leq l\neq r\leq q-1.
		\end{eqnarray*}
\end{proposition}
\begin{proof}
		It can be easily proved. Hence, we omit it. 
	\end{proof}
	We know that $\{\phi^{(j)}_{k,l}: k\in \N_{0},0\leq l\leq q-1, j\in {\Z}\}$ is an orthonormal system if and only if $M^{(j)}(\xi)$ is a unitary matrix. At first, we find equivalent conditions for $M^{(j)}(\xi)$ to be unitary. It is analogous to see that
	$M^{(j)}(\xi)$ is unitary if and only if
	\begin{eqnarray*}
		\begin{aligned}
			(a) \sum_{\mu=0}^{q-1} m^{(j)}_{l}|(\mathfrak{p}\xi+\mathfrak{p}\lambda(\mu))|^2=1 \ \text{a.e. for all} \ l=0,1,2,...,q-1.\qquad\qquad\qquad\quad\quad \ \label{filter0}\\
			(b)\sum_{\mu=0}^{q-1} m^{(j)}_{l}|(\mathfrak{p}\xi+\mathfrak{p}\lambda(\mu))|^2  \overline{m^{(j)}_{r}|(\mathfrak{p}\xi+\mathfrak{p}\lambda(\mu))|^2}=0 \ \ \text{a.e. for all} \ 0\leq l\neq r\leq q-1.\\
		\end{aligned} \label{eq:2.4}
	\end{eqnarray*}
 and scale relations, including the density condition, provide
	\begin{equation*}
	\lim_{j \rightarrow +\infty }|m_0^{(j)}(\mathfrak{p}^{j} \xi)|=1.
	\end{equation*}

	\begin{theorem}
		Let $\{ {\rm e}_t: t\in {\N_{0}} \}$ be an orthonormal system in $H^{\mathfrak s}(K_q)$. Then the sequence $\{\phi^{(j)}_{k,l}: k\in \N_{0},0\leq l\leq q-1, j\in {\Z}\}$ given by
		\eqref{eq:2.1} is an orthonormal system in $H^{\mathfrak s}(K_q)$ if and only if the matrix $M^{(j)}(\xi)$ given by \eqref{eq:2.3} is unitary for
		a.e. $\xi \in K_q$.
	\end{theorem}
	\begin{proof}
		Using definition \ref{def2}, we get
		\begin{eqnarray*}
			&&\bigl\langle \phi^{(j)}_{k,l},\phi^{(j)}_{m,l}\bigr\rangle_{H^{\mathfrak s}(K_q)}= \int_{K_q} (1+|\xi|^2)^{\mathfrak s} \hat{\phi}^{(j)}_{k,l}(\xi) \overline{\hat{\phi}^{(j)}_{m,l}(\xi)} {\rm d}\xi\\
			&&=q \int_{K_q} (1+|\xi|^2)^{\mathfrak s}  \sum_{ t\in {\N_{0}}} \alpha^{(j)}_{2k-t,l} \overline{\alpha^{(j)}_{2m-t,l}} {\hat{{\rm e}_{t,l}}}(\xi) \overline{\hat{{\rm e}_{t,l}}(\xi)} {\rm d}\xi\\
			&&=q \sum_{t\in {\N_{0}}} \alpha^{(j)}_{2k-t,l} \overline{\alpha^{(j)}_{2m-t,l}} \bigl\langle{{{\rm e}_{t,l}}}(\xi) , {{{\rm e}_{t,l}}}(\xi)\bigr\rangle_{\mathfrak s}\\
			&&=q \sum_{t\in {\N_{0}}} \alpha^{(j)}_{2k-t,l} \overline{\alpha^{(j)}_{2m-t,l}} \\
			&&=q \sum_{t\in {\N_{0}}} \alpha^{(j)}_{\beta-2(m-k),l} \overline{\alpha^{(j)}_{\beta}} \\
			&&=\delta_{2(m-k),0}^{l}=\delta_{2m,2k}^{l}=\delta^{l}_{m,k},
		\end{eqnarray*}
		for all $0\leq k,m\leq q-1$. Further, we have
		\begin{eqnarray*}
			&&\bigl\langle \phi^{(j)}_{k,l},\phi^{(j)}_{m,r}\bigr\rangle_{H^{\mathfrak s}(K_q)}= \int_{K_q} (1+|\xi|^2)^{\mathfrak s} \hat{\phi}^{(j)}_{k,l}(\xi) \overline{\hat{\phi}^{(j)}_{m,r}(\xi)} {\rm d}\xi\\
			&&=q \int_{K_q} (1+|\xi|^2)^{\mathfrak s}  \sum_{ t\in {\N_{0}}} \alpha^{(j)}_{2k-t,l} \overline{\alpha^{(j)}_{2m-t,r}} {\hat{{\rm e}_{t,l}}}(\xi) \overline{\hat{{\rm e}_{t,r}}(\xi)} {\rm d}\xi\\
			&&=q \sum_{t\in {\N_{0}}} \alpha^{(j)}_{2k-t,l} \overline{\alpha^{(j)}_{2m-t,r}} \bigl\langle{{{\rm e}_{t,l}}}(\xi) , {{{\rm e}_{t,r}}}(\xi)\bigr\rangle_{\mathfrak s}\\
			&&=q \sum_{t\in {\N_{0}}} \alpha^{(j)}_{2k-t,l} \overline{\alpha^{(j)}_{2m-t,r}} \\
			&&=q \sum_{t\in {\N_{0}}} \alpha^{(j)}_{\beta-2(m-k),l} \overline{\alpha^{(j)}_{\beta,r}}=0.
		\end{eqnarray*}
		This proves that $\{\phi^{(j)}_{k,l}: k\in \N_{0},0\leq l\leq q-1, j\in {\Z}\}$ is an orthonormal system. The converse part of this theorem can be proved by reversing the above steps.
	\end{proof}
	
	\begin{lemma}\label{lem:split}
		Let $\varphi^{(j)} \in H^{\mathfrak s}(K_q)$ be such that $\{\varphi^{(j)}(\cdot-\lambda(k)): k\in \mathbb{N}_0\}$ is an orthonormal system.
		Let $\mathsf{V} = \overline{\rm span}\{q^{1/2}\varphi^{(j)}(\mathfrak{p}^{-1}\cdot-\lambda(k)): k\in \mathbb{N}_0\}$. Let
		$m^{(j)}_l(\xi) = q^{-1/2}\sum _{k=0}^{\infty}\alpha^{(j)}_{k,l}  \overline{\chi_k(\xi)}$, $0\leq l\leq q-1$, where
		$\{\alpha^{(j)}_{k,l}: k\in \mathbb{N}_0\}\in\ell^2(\mathbb{N}_0)$ for $0\leq l\leq q-1$. Define ${\hat\psi}^{(j)} _{l}(\xi) = m^{(j+1)}_l(\mathfrak{p}\xi)\hat{\varphi}^{(j+1)}(\mathfrak{p}\xi)$. Then $\{\psi^{(j)}_l (\cdot - \lambda(k)): 0\leq l\leq q-1, k\in \N_0\}$ is an orthonormal system in $\mathsf{V}$ if and only if the matrix
		\begin{equation*}
		M^{(j+1)}(\xi) = \Bigl(m^{(j+1)}_l(\mathfrak{p}\xi + \mathfrak{p}\lambda(k))\Bigr)_{l,k=0}^{q-1}
		\end{equation*}
		is unitary for a.e. $\xi\in\mathcal{D}_q$.
	\end{lemma}
	
	\proof Suppose $M^{(j+1)}(\xi)$ is unitary for a.e $\xi\in\mathcal{D}_q$. Then, for $0\leq r, t\leq q-1$ and $k, l\in\N_0$, we get
	\begin{eqnarray*}
		\lefteqn{\Bigl\langle\psi^{(j)}_r\bigl(\cdot - \lambda(k)\bigr),\psi^{(j)}_t\bigl(\cdot - \lambda(l)\bigr)\Bigr\rangle_{H^{\mathfrak s}(K_q)}} \\
		& = & \int_{\mathcal{D}_q}\sum\limits_{\mu=0}^{q-1}\sum\limits_{n\in\N_0}(1+\mathfrak{p}^{-2}(\mathfrak{p}|\xi|+\mathfrak{p}\lambda(qn+\mu))^2)^{\mathfrak s}m^{(j+1)}_r
		\bigl(\mathfrak{p}\xi+\mathfrak{p}\lambda(qn+\mu)\bigr)
		\overline{m^{(j+1)}_t\bigl(\mathfrak{p}\xi+\mathfrak{p}\lambda(qn+\mu)\bigr)} \\
		&   & \qquad \times\bigl|\hat\varphi^{(j+1)}\bigl(\mathfrak{p}\xi+\mathfrak{p}\lambda(qn+\mu)\bigr)\bigr|^2
		\overline{\chi_k(\xi)}\chi_l(\xi)~{\rm d}\xi \\
		& = & \delta_{r,t}\delta_{k,l}.
	\end{eqnarray*}
	Hence, $\{\psi^{(j)}_l(\cdot-\lambda(k)): 0\leq l \leq q-1, k\in\N_0\}$ is an orthonormal system in $\mathsf{V}$. Similarly, the converse of this Lemma can be proved by reversing the above steps.
	
	Now we simplify the index. For $n\in \N_{0}$, we can write the expansion of $n$ uniquely in the form of base $q$, we have unique $\mu=(\mu_1,\mu_2,\dots,\mu_j)$ such that 
	\begin{equation}\label{eqn:q.ary}
	n=\mu_1+\mu_2q+\mu_3q^2+\cdots+\mu_jq^{j-1}, \ \text{for} \ j=1,2,3,\dots,
	\end{equation}
	where $0\leq\mu_i\leq q-1$ for all $i=1,2,\dots ,j$ and $\mu_j\not=0$.
	Observe that \eqref{eqn:q.ary} is always a finite sum. Indeed, if $q^{d_0-1}\leq n < q^{d_0}$, say, then we have $n=q^{d_0-1}+ n_{1}$  where $q^{d_1-1}\leq n_1 < q^{d_1}$ and $d_1<d_0$; and iterating this procedure, we obtain $n=q^{d_0-1}+ q^{d_1+1}+\dots+q^{d_k+1}$, where $1\leq d_k<d_{k-1}<\dots <d_0$. That is, $\mu_{j}=1$ for $j=d_k,\dots, d_0$, and $\mu_{j}=0$, otherwise. 
	 When we order the elements of the set like $\{0,1,2,\dots,q-1\}$ in any way, we can write $\mu_i \in \{0,1,2,\dots,q-1\}, \forall i$. Let $\Delta_j$ be the set of these $j$-tuple $\mu=(\mu_1,\mu_2,\dots,\mu_j)$ with length $j$, and denote $\Delta=\bigcup_{j=1}^{+\infty} \Delta_j$. Notice that when $i\leq j, \Delta_i$ can be embedded in $\Delta_j$ naturally, by considering $(\mu_1,\mu_2,\dots,\mu_j)$ as $(\mu_1,\mu_2,\dots,\mu_j,0,\dots,0)$.
	
	Suppose $w^{(j)}_0=\varphi^{(j)}$. Now we can define the basic wavelet packets associated with the scaling function $\varphi^{(j)}$ recursively as follows: \\
	Suppose $w_{\ell}$ is defined for $\ell\geq 0$ and for $0 \leq r \leq q-1$, define
	\begin{equation}\label{wpkt}
	{w}^{(j)}_n(x)=w^{(j)}_{r+q\ell}(x) = q\sum\limits_{k\in \N_0}\alpha^{(j+1)}_{k,r} w^{(j+1)}_{\ell}(\mathfrak{p}^{-1}x-\lambda(k)). \end{equation}
	Note that equation \eqref{wpkt} defines $w_n$ for every integer $n\geq 0$. Applying Fourier transform on equation \eqref{wpkt}, we obtain
	\begin{equation}\label{e.ftwpkt}
	\hat{w}^{(j)}_n(\xi)=(w^{(j)}_{r+q\ell})^\wedge (\xi) = m^{(j+1)}_{\mu_r}(\mathfrak{p}\xi)\hat{w}^{(j+1)}_{\ell}(\mathfrak{p}\xi).
	\end{equation}
	
	\begin{proposition} \label{sec:pro:4}
		For $n\in \N_{0}$ and let the $q$-adic expansion of $n$ be given by \eqref{eqn:q.ary}. Then the Fourier transform of the basic wavelet
		packet defined by \eqref{e.ftwpkt} is given by
		\begin{eqnarray}
		\hat{w}_{n}^{(j)} (\xi)= \prod_{J=1}^{+\infty} m_{\mu_{r}}^{(j+J)} (\mathfrak{p}^{J} \xi) (1+\mathfrak{p}^{-2j} |\xi|^2)^{-s/2}, \quad \xi\in K_q. \label{eq:2.8}
		\end{eqnarray}
	\end{proposition}
	
	\begin{proof}
		The Fourier transform equivalence of the $q$-adic scale relations \eqref{wpkt} for the wavelet packets is given by \eqref{e.ftwpkt}.
		We prove \eqref{eq:2.8} by induction on $n:=q\ell+r$ where $\ell\in {\N_{0}}$ and $r=0,1,2,\dots,q-1$. Suppose \eqref{eq:2.8} holds for all $n$ with
		$0\leq n< q^{d_0}$ and consider $q^{d_0}\leq n < q^{d_0+1}$ with the $q$-adic expansion \eqref{e.ftwpkt}. From the discussion presented before, we have
		\begin{eqnarray*}
			\mu_r=
			\begin{cases}
				{1,2,\dots,q-1} & \mbox{for} \ r=d_0+1\\
				0 & \mbox{for} \ r>d_0+1.
			\end{cases}
		\end{eqnarray*}
		So that
		\begin{equation*}
		n=\sum_{r=1}^{d_0+1} \mu_r q^{r-1},
		\end{equation*}
		and
		\begin{equation*}
		{\frac{n}{q}}=\frac{\mu_{1}}{q}+\sum_{r=1}^{d_0} \mu_{r+1} q^{r-1}.
		\end{equation*}
		Let the symbol $\floor*{x}$ denote the greatest integer not exceeding $x$. Notice that
		\begin{equation*}
		n=q\floor*{\frac{1}{q}} + \mu_{1}.
		\end{equation*}
		Hence, from \eqref{wpkt}, we have
		\begin{eqnarray}
		\hat{w}^{(j)}_n(\xi)=m^{(j+1)}_{\mu_{1}} (\mathfrak{p}\xi) \hat{w}^{(j+1)}_{\floor*{\frac{n}{q}}}(\mathfrak{p}\xi). \label{eq:2.10}
		\end{eqnarray}
		On the right-hand-side of equation \eqref{eq:2.10}, we can write $\floor*{\frac{n}{q}}$ as 
		\begin{equation*}
		\floor*{\frac{n}{q}}=\sum_{r=1}^{d_0} \mu_{r+1} q^{r-1}\leq \frac{n}{q}<q^{d_0},
		\end{equation*}
		it follows the induction hypothesis that
		\begin{eqnarray}
		\hat{w}^{(j+1)}_{\floor*{\frac{n}{q}}}(\xi)= m^{(j+2)}_{\mu_2} (\mathfrak{p}^{2}\xi) \hat{w}^{(j+2)}_{\floor*{\frac{n}{q^2}}}(\mathfrak{p}^2\xi). \label{eq:2.11}
		\end{eqnarray}
		Then from \eqref{eq:2.10} and \eqref{eq:2.11}, we have
		\begin{eqnarray*}
			&&\hat{w}^{(j)}_n(\xi)=m^{(j+1)}_{\mu_1} (\mathfrak{p}\xi) m^{(j+2)}_{\mu_2} (\mathfrak{p}^2\xi) \hat{w}^{(j+2)}_{\floor*{\frac{n}{q^2}}}(\mathfrak{p}^2\xi)\\
			&&=\prod_{J=1}^{2} m^{(j+J)}_{\mu_r} (\mathfrak{p}^J\xi)  \hat{w}^{(j+2)}_{\floor*{\frac{n}{q^2}}}(\mathfrak{p}^{2}\xi)\\
			&&=\dots=\prod_{J=1}^{k} m^{(j+J)}_{\mu_r} (\mathfrak{p}^J\xi)  \hat{w}^{(j+k)}_{\floor*{\frac{n}{q^k}}}(\mathfrak{p}^{k}\xi)\\
				\end{eqnarray*}
			Using the density condition provided in Section \ref{sec:3.2} on the right-hand side of the above expression, we have the desired result as follows
		\begin{eqnarray*}
			\hat{w}^{(j)}_n(\xi)= \prod_{J=1}^{+\infty} m_{\mu_r}^{(j+J)} (\mathfrak{p}^{J} \xi)(1+\mathfrak{p}^{2j} |\xi|^2)^{-{\mathfrak s}/2}.
		\end{eqnarray*}
	
	\end{proof}
	
	\section{Wavelet packet functions generated by MRA in $H^{\mathfrak s}( K_q)$ at different scales}
If $\varphi^{(j)}$ generates an orthonormal MRA $\{\mathsf{V}_j\}_{j\in {\mathbb{Z}}}$ associated with wavelet function
	$\psi^{(j)}$ in $H^{\mathfrak s}(K_q)$. Then we define the wavelet packets as sequence $\{w_n^{(j)}: n\geq0\}$ of functions
	$w_0^{(j)}=\varphi^{(j)}, and \ w_n^{(j)}=\psi_n^{(j)} \ (1\leq n\leq q-1)$, where $\hat{\psi}_{\ell}(\xi)=m(\mathfrak{p}\xi) \hat{\varphi}^{(j)}(\mathfrak{p}\xi)\ (1\leq n\leq q-1)$. Suppose $w_{\ell}$ is defined for ${\ell}\geq 0$ and $0\leq r \leq q-1$, we have
	\begin{eqnarray}
	\begin{aligned}
	w^{(j)}_{n}(\cdot)=w^{(j)}_{r+q\ell}(\cdot)=q^{j/2} \sum_{r\in \N_{0}} \alpha^{(j+1)}_{k,r} w^{(j+1)}_{\ell} (\mathfrak{p}^{-j}\cdot-\lambda(k)).
	\end{aligned} \label{eq:3.1}
	\end{eqnarray}
	\begin{definition} Suppose $n\in {\mathbb{N}_0}$ and let $w^{(j)}_n$ be a wavelet function involving scaling function $\varphi^{(j)}$. Then, for any integer
		$j$ and $k\in \N_{0}$, we define
		\begin{eqnarray}
		w^{(j)}_{j,k,n}(\cdot)=q^{j/2} w^{(j)}_n(\mathfrak{p}^{-j}\cdot-\lambda(k)). \label{eq:3.2}
		\end{eqnarray}
	\end{definition}
	
	\begin{lemma}
		If ${\mathfrak s} \in {\mathbb{R}}$,  $w^{(j+1)}_{j,k,n}\in H^{\mathfrak s}(K_q)$, $j\in \Z$ and $k,n\in \N_{0}$, then the distributions
		$ \left\{q^{\frac{(j)}{2}} w^{(j+1)}_{\floor*{\frac{n}{q}}}(\mathfrak{p}^{(j+1)}\cdot-\lambda(k))\right\}$, are orthonormal in $H^{\mathfrak s}(K_q)$ if and only if
		
		\begin{eqnarray*}
			\sum\limits_{l_1\in\N_0}(1+\mathfrak{p}^{-2(j+1)}(\mathfrak{p}|\xi|+\mathfrak{p}\lambda(l_1))^2)^{\mathfrak s} \hat{w}_{\floor*{\frac{n}{q}}}^{(j+1)}(\mathfrak{p}\xi+\mathfrak{p}\lambda(l_1))\overline{\hat{w}_{\floor*{\frac{m}{q}}}^{(j+1)}(\mathfrak{p}\xi+\mathfrak{p}\lambda(l_1))}=\delta_{{\floor*{\frac{n}{q}}},{\floor*{\frac{m}{q}}}}.
		\end{eqnarray*}
		
	\end{lemma}
	\begin{proof}
		Since $w^{(j+1)}_{j,k,n}\in H^{\mathfrak s}(K_q)$, the series
		\begin{eqnarray*}
			S^{(j+1)}(\mathfrak{p}\xi)=\sum\limits_{l_1\in\N_0}(1+\mathfrak{p}^{-2(j+1)}(\mathfrak{p}|\xi|+\mathfrak{p}\lambda(l_1))^2)^{\mathfrak s} \hat{w}_{\floor*{\frac{n}{q}}}^{(j+1)}(\mathfrak{p}\xi+\mathfrak{p}\lambda(l_1))\overline{\hat{w}_{\floor*{\frac{m}{q}}}^{(j+1)}(\mathfrak{p}\xi+\mathfrak{p}\lambda(l_1))},
		\end{eqnarray*}
		converges a.e., belongs to $L^1(\mathcal{D}_q)$, and is integral-periodic. Further, for every $l\in {\N_{0}}$, we have
		\begin{align*}
			\int_{\mathcal{D}_q} S^{(j+1)}(\mathfrak{p}\xi) \overline{\chi_{k-l}(\mathfrak{p}\xi)} {\rm d}\xi&=\int_{\mathcal{D}_q}\sum\limits_{l_1\in\N_0}(1+\mathfrak{p}^{-2(j+1)}(\mathfrak{p}|\xi|+\mathfrak{p}\lambda(l_1))^2)^{\mathfrak s} \hat{w}_{\floor*{\frac{n}{q}}}^{(j+1)}(\mathfrak{p}\xi+\mathfrak{p}\lambda(l_1))\\ 
			&\quad\times\overline{\hat{w}_{\floor*{\frac{m}{q}}}^{(j+1)}(\mathfrak{p}\xi+\mathfrak{p}\lambda(l_1))} \overline{\chi_{k-l}(\mathfrak{p}\xi)} {\rm d}\xi\\
			&=\int_{K_q}(1+\mathfrak{p}^{-2j}|\xi|)^2)^{\mathfrak s} \hat{w}_{\floor*{\frac{n}{q}}}^{(j+1)}(\mathfrak{p}\xi) \overline{\hat{w}_{\floor*{\frac{m}{q}}}^{(j+1)}(\mathfrak{p}\xi)} \overline{\chi_{k-l}(\mathfrak{p}\xi)} {\rm d}\xi\\
			&=\Bigl\langle q^{-j/2}{w}_{\floor*{\frac{n}{q}}}^{(j+1)}(\mathfrak{p}^{(j+1)}x-\lambda(l)),q^{-j/2}{w}_{\floor*{\frac{m}{q}}}^{(j+1)}(\mathfrak{p}^{(j+1)}x-\lambda(k))\Bigr\rangle_{H^{\mathfrak s}(K_q)} \\
			&=\delta_{{\floor*{\frac{n}{q}}},{\floor*{\frac{m}{q}}}}.\\
		\end{align*}
		Therefore,
		\begin{eqnarray*}
			\int_{\mathcal{D}_q} S^{(j+1)}(\mathfrak{p}\xi) \overline{\chi_{k-l}(\mathfrak{p}\xi)} d\xi=\delta_{{\floor*{\frac{n}{q}}},{\floor*{\frac{m}{q}}}},
		\end{eqnarray*}
		if and only if $S^{(j+1)}(\mathfrak{p}\xi)=\delta_{{\floor*{\frac{n}{q}}},{\floor*{\frac{m}{q}}}}$ and $k=l$.
	\end{proof}
	
	\begin{theorem}[Orthonormality at $j^{th}$ level] \label{sec:them:3}
		Let $j$ be an integer and for all $k,l,m$ and $n$ belongs to $\N_{0}$, then we have
		\begin{equation}
		\Bigl\langle w^{(j)}_{j,k,n}\bigl(\cdot\bigr),w^{(j)}_{j,l,m}\bigl(\cdot\bigr)\Bigr\rangle_{H^{\mathfrak s}(K_q)}=\delta_{{\floor*{\frac{n}{q}}},{\floor*{\frac{m}{q}}}}\delta_{\mu_{1},\eta_{1}}\delta_{k,l},  \label{eq:3.5}
		\end{equation}
		where $n=\mu_{1}+\mu_{2}q+\dots+\mu_{j}q^{j-1}$, $1\leq \mu_{i}\leq q-1$ for $i=1,2,3,\dots, j$ with $\mu_{j}\neq 0$ and $m=\eta_{1}+\eta_{2}q+\dots+\eta_{j} q^{j-1}$, $1\leq \eta_{i}\leq q-1$ for $i=1,2,3,\dots, j$ with $\eta_{j}\neq 0$.
	\end{theorem}
	
	\begin{proof}
		For all $k,l,m$ and $n$ belongs to $\N_{0}$, we have
		\begin{align*}
			\lefteqn{\Bigl\langle w^{(j)}_{j,k,n}\bigl(\cdot\bigr),w^{(j)}_{j,l,m}\bigl(\cdot\bigr)\Bigr\rangle_{H^{\mathfrak s}(K_q)}}\\
			& =  \int_{\mathcal{D}_q}\sum\limits_{\mu=0}^{q-1}\sum\limits_{l_1\in\N_0}(1+\mathfrak{p}^{-2(j+1)}(\mathfrak{p}|\xi|+\mathfrak{p}\lambda(ql_1+\mu))^2)^{\mathfrak s} m^{(j+1)}_{\mu_{1}}
			\bigl(\mathfrak{p}\xi+\mathfrak{p}\lambda(ql_1+\mu)\bigr)\\ 
			&\quad \times \hat{w}_{\floor*{\frac{n}{q}}}^{(j+1)}(\mathfrak{p}\xi+\mathfrak{p}\lambda(ql_1+\mu))\overline{\chi_k(\xi)}\overline{m^{(j+1)}_{\eta_{1}}\bigl(\mathfrak{p}\xi+\mathfrak{p}\lambda(ql_1+\mu)\bigr)}\\
			&\quad \times    \overline{\hat{w}_{\floor*{\frac{m}{q}}}^{(j+1)}(\mathfrak{p}\xi+\mathfrak{p}\lambda(ql_1+\mu))}
			\chi_l(\xi)~{\rm d}\xi \\
			& =  \delta_{{\floor*{\frac{n}{q}}},{\floor*{\frac{m}{q}}}}\delta_{\mu_{1},\eta_{1}}\delta_{k,l}.
		\end{align*}
		Hence, $\{w^{(j)}_{n,k,l}(\cdot): 0\leq n \leq q-1, l, k\in\N_0\}$ is an orthonormal system in $\mathsf{V}$. Likewise, the converse of this theorem can be proved by reversing the above steps.
	\end{proof}
	The following proposed theorem provides an alternate form for orthogonality of wavelet packets at $j^{th}$ level in $H^{\mathfrak s}(K_q)$ by exploiting the theory of convolution of Fourier transforms.   
	
	\begin{theorem}
		 The Fourier transform of $\kappa_{-s/2}$ is given by  $\hat{\kappa}_{-s/2}(\xi)=(1+|\xi|^2)^{-{{\mathfrak s}/2}}$ (i.e., replace $\theta=\mathfrak{-s/2}\in \mathbb{R}$ in example \ref{ex7}) and $w^{(j)}_{n,k}(\cdot)=q^{j/2}w^{(j)}_{n}(\mathfrak{p}^{-j}\cdot-\lambda(k))$ where $j\in {\Z}$ and $n,k\in \N_{0}$. Then
		\begin{eqnarray*}
			\Bigl\langle\kappa_{-s/2}\ast w^{(j)}_{n,k},\kappa_{-s/2}\ast w^{(j)}_{m,l}\Bigr\rangle_{H^{\mathfrak s}(K_q)}=\delta_{{\floor*{\frac{n}{q}}},{\floor*{\frac{m}{q}}}}\delta_{\mu_{1},\eta_{1}}\delta_{k,l}.
		\end{eqnarray*}
		
	\end{theorem}
	\begin{proof} We use the convolution theorem for Fourier transforms and obtain the following
		\begin{align*}
			&\Bigl\langle\kappa_{-s/2} \ast w^{(j)}_{n,k},\kappa_{-s/2} \ast w^{(j)}_{m,l}\Bigr\rangle_{\mathfrak s}
			&= \int_{K_q} (1+|\xi|^2)^{\mathfrak s} {\hat\kappa_{-s/2}}(\xi) {{\hat w}^{(j)}_{n,k}}(\xi) \overline{{\hat\kappa_{-s/2}}(\xi)} \ \overline{{{\hat w}^{(j)}_{m,l}}(\xi)} {\rm d}\xi.
		\end{align*}
		Since ${\hat\kappa_{-s/2}}(\xi)=(1+|\xi|^2)^{-{\mathfrak s}/2}$, then the right-hand side of the above expression becomes
		\begin{equation*}
			\int_{K_q}  {{\hat w}^{(j)}_{n,k}}(\xi)  \overline{{{\hat w}^{(j)}_{m,l}}(\xi)} d\xi.
		\end{equation*}
		By Parseval formula and relation \eqref{eq:3.2}, we get
		\begin{eqnarray*}
			\Bigl\langle\kappa_{-s/2} \ast w^{(j)}_{n,k},\kappa_{-s/2} \ast w^{(j)}_{m,l}\Bigr\rangle_{H^{\mathfrak s}(K_q)} 
			&=&\delta_{{\floor*{\frac{n}{q}}},{\floor*{\frac{m}{q}}}}\delta_{\mu_{1},\eta_{1}}\delta_{k,l}.
		\end{eqnarray*}
	\end{proof}
We are interested in constructing exciting examples of wavelet packets in $H^{\mathfrak s}(K_q)$ and its subsets. We provide three compelling examples: Haar wavelet packets, Weierstrass wavelet packets, and 3-adic Cantor wavelet packets. The Weierstrass and 3-adic Cantor wavelet packets are from a family of fractal wavelet packets.
	\begin{example}\label{ex8}
		For all $k,l,m$ and $n$ belongs to $\N_{0}$ and the wavelet packet function defined by 
		\begin{eqnarray}\label{eqe8}
			w^{(j)}_{j,k,n}(x)=\kappa_{-s/2}\bigl(x-\lambda(k)\bigr)\ast 	q^{j}\Phi^{(j)}_{\mathcal{P}_{q}^{j},n}{\bigl(x-\lambda(k)\bigr)}.
		\end{eqnarray} 
	Moreover, this wavelet packet functions is also orthogonal, i.e., 
	\begin{align*}
	\Bigl\langle w^{(j)}_{j,k,n}\bigl(\cdot\bigr),w^{(j)}_{j,l,m}\bigl(\cdot\bigr)\Bigr\rangle_{H^{\mathfrak{s}}(K_q)}=\delta_{k,l}.
	\end{align*}
\end{example}
\begin{solution}
	By Theorem \ref{sec:them:3}, for all $k,l,m$ and $n$ belongs to $\N_{0}$, we have
	\begin{align*}
		\lefteqn{\Bigl\langle w^{(j)}_{j,k,n}\bigl(\cdot\bigr),w^{(j)}_{j,l,m}\bigl(\cdot\bigr)\Bigr\rangle_{H^{\mathfrak{s}}(K_q)}}\\
			& =  \int_{K_q} (1+|\xi|^2)^{\mathfrak s}(1+|\xi|^2)^{\mathfrak -s} \overline{\chi_k(\xi)}\chi_l(\xi){\Phi}^{(j)}_{{\Gamma}^{j},n}{\bigl(\xi\bigr)}{\Phi}^{(j)}_{{\Gamma}^{j},m}{\bigl(\xi\bigr)}~{\rm d}\xi\\
				&=\delta_{k,l}
	\end{align*}

	\end{solution}
	\begin{example}\label{ex9}[Weiestrass type wavelet packets]
		For all $k,l,m$ and $n$ belongs to $\N_{0}$ and the wavelet packets on $H^{\mathfrak s}( K_2)$ defined by 
		\begin{eqnarray}\label{eqe9}
			w^{(j)}_{j,k,n}(x)=\kappa_{-s/2}\bigl(x-\lambda(k)\bigr)\ast 2^2 {\mathcal W}^{(j)}_{2,n}(x-\lambda(k)).
		\end{eqnarray} 
		Moreover, these wavelet packets are also orthogonal, i.e., 
		\begin{align*}
			\Bigl\langle w^{(j)}_{j,k,n}\bigl(\cdot\bigr),w^{(j)}_{j,l,m}\bigl(\cdot\bigr)\Bigr\rangle_{H^{\mathfrak s}(K_2)}=\delta_{k,l}.
		\end{align*}
	\end{example}
	\begin{solution}
		By Theorem \ref{sec:them:3}, for all $k,l,m$ and $n$ belongs to $\N_{0}$, we have
		\begin{align*}
			\lefteqn{\Bigl\langle w^{(j)}_{j,k,n}\bigl(\cdot\bigr),w^{(j)}_{j,l,m}\bigl(\cdot\bigr)\Bigr\rangle_{H^{\mathfrak s}(K_2)}}\\
				& = 2^4 \int_{{\mathcal D}_2} (1+|\xi|^2)^{\mathfrak s} \hat{\kappa}_{{-s/2}}^{(j)}(\xi)\hat{{\mathcal W}}^{(j)}_{2,n}\bigl(\sum_{j=s}^{+\infty}x_j{\mathfrak{2^{-j}}}-\lambda(k)\bigr)(\xi)\overline{\hat{\kappa}_{-s/2}^{(j)}(\xi)}~\overline{\hat{{\mathcal W}}^{(j)}_{2,n}\bigl(\sum_{j=s}^{+\infty}x_j{\mathfrak{2^{-j}}}-\lambda(l)\bigr)(\xi)}{\rm d}\xi\\
		\end{align*}
	by using a change of the variable technique on the right-hand side of the above expression, we have
	\begin{align*}
			& =  2^4 \int_{{\mathcal D}_2} (1+|\xi|^2)^{\mathfrak s} |\hat{\kappa}_{{-s/2}}^{(j)}(\xi)|^2\hat{{\mathcal W}}^{(j)}_{2,n}\bigl(\xi\bigr)\overline{\chi_k(\xi)}~\overline{\hat{{\mathcal W}}^{(j)}_{2,n}\bigl(\xi\bigr)}\chi_k(\xi){\rm d}\xi
		\end{align*}
	Using equations \eqref{eqe5.1} and \eqref{eqe7.1} on the above expression, we have 
	\begin{align*}
				& =  \int_{{\mathcal D}_2} (1+|\xi|^2)^{\mathfrak s} (1+|\xi|^2)^{\mathfrak -s}\overline{\chi_k(\xi)}\chi_l(\xi)~{\rm d}\xi\\
			&=\delta_{k,l}.
		\end{align*}
	\end{solution}

\begin{example}\label{ex10}[3-adic Cantor type wavelet packets]
	For all $k,l,m$ and $n$ belongs to $\N_{0}$ and the wavelet packets on $H^{\mathfrak s}(K_3)$ defined by 
	\begin{eqnarray}\label{eqe10}
		w^{(j)}_{j,k,n}(x)=\kappa_{-s/2}\bigl(x-\lambda(k)\bigr)\ast 2 {\mathcal C}^{(j)}_{n}(x-\lambda(k)).
	\end{eqnarray} 
	Moreover, these wavelet packets are also orthogonal, i.e., 
	\begin{align*}
		\Bigl\langle w^{(j)}_{j,k,n}\bigl(\cdot\bigr),w^{(j)}_{j,l,m}\bigl(\cdot\bigr)\Bigr\rangle_{H^{\mathfrak s}(K_3)}=\delta_{k,l}.
	\end{align*}
\end{example}
\begin{solution}
	One can prove this wavelet satisfies an orthogonality condition similar to the example.
\end{solution}
\section{Conclusion}
 In this work, we define Sobolev spaces $H^{\mathfrak{s}}(K_q)$ over a local field $K_q$ of finite characteristic $p>0$. We explore the application of fractal functions, specifically the Weierstrass and 3-adic Cantor type, to illustrate the potential of these spaces. A Multi-Resolution Analysis (MRA) is developed, and wavelet expansions are investigated, focusing on the orthogonality of wavelet packet functions at various scales. Using convolution theory and Fourier transforms, we construct Haar and fractal wavelet packets within $H^{\mathfrak{s}}(K_q)$, proving their orthogonality via Theorem \ref{sec:them:3}. The findings underscore the utility of wavelet packets in mathematical science and engineering.

\section*{Acknowledgments} The author gratefully acknowledges Prof. Peter Robert Massopust, Privatdozent Center of Mathematics at the Technical University of Munich, Germany, for his contributions. Prof. Massopust's insightful conversations and expert direction shaped the practical demonstration of  Sobolev space through examples. Prof. Massopust painstakingly reviewed the manuscript, demonstrating his dedication and scholarly zeal. His perceptive feedback and scholarly commitment have improved the 
clarity and sophistication of this research work academic.

\section*{Declaration of interests} The authors declare that they have no known competing financial interests or personal relationships that could have appeared to influence the work reported in this paper.
		\bibliographystyle{plain}
	\bibliography{Ref_Local}

\end{document}